\documentclass[12pt]{article} 
    \usepackage{graphicx} 

    \usepackage{epsfig} 
    \usepackage{amsmath} 
    \usepackage{amssymb} 
    \usepackage{amsthm} 
    \usepackage{latexsym} 
    \usepackage{cite}

    \newcommand{\comment}[1]{}

\begin{document} 
      
  \title{\LARGE 
  {\bf Continuum Nonsimple Loops \\ and $2$D Critical Percolation} 
  } 

  \author{ 
  {\bf Federico Camia} 
  \thanks{EURANDOM, P.O. Box 513, 5600 MB Eindhoven, The Netherlands}\,
  \thanks{E-mail: camia@eurandom.tue.nl}\\
  \and 
  {\bf Charles M.~Newman} 
  \thanks{Courant Inst.~of Mathematical Sciences, 
  New York University, New York, NY 10012, USA}\,
  \thanks{E-mail: newman@courant.nyu.edu}\,
  \thanks{Research partially supported by the 
  U.S. NSF under grant DMS-01-04278.}
  } 

  \date{} 

  \maketitle

\begin{abstract} 
Substantial progress has been made in recent years on the 2D critical 
percolation 
scaling limit and its conformal invariance properties. 
In particular, chordal $SLE_6$ (the Stochastic L\"{o}wner 
Evolution with parameter $\kappa = 6$) was, in the work of
Schramm and of Smirnov, identified as
the scaling limit of the critical percolation ``exploration process.'' 
In this paper we use that and other results to construct what we argue 
is the \emph{full} scaling limit of the collection of \emph{all} closed 
contours surrounding the critical percolation clusters on the 2D triangular 
lattice. 
This random process or gas of continuum nonsimple loops in 
${\mathbb R}^2$ is constructed 
inductively by repeated use of chordal $SLE_6$. 
These loops do not cross but do touch each other --- 
indeed, any two loops are connected by a finite ``path'' of touching loops. 
\end{abstract} 

\noindent {\bf Keywords:} scaling limit, percolation, SLE, 
continuum loops, nonsimple loops, triangular lattice, conformal invariance. 


\section{Introduction} \label{intro}


Percolation is a model with a wide range of applications and, 
especially in two dimensions, a well developed theory (see, for example, 
 \cite{kesten,grimmett}). 
It has been used as a proving ground for developing tools that
can be applied 
to more complex systems, and is of great interest in its own right, as 
it is 
perhaps the simplest (non-mean-field) 
model displaying a phase transition with features 
such as scaling and universality at criticality. 

    In the critical case, the fractal and conformally invariant nature of 
  (the scaling limit of) large 
    percolation clusters has attracted much attention and is of interest 
  for 
  both 
    intrinsic reasons and as a paradigm for the study of other systems. 

    The ground-breaking work of Schramm~\cite{schramm} 
    and Smirnov~\cite{smirnov} 
    has elucidated much about the nature of the 
    scaling limit of the cluster 
    boundaries or contours in terms of $SLE_6$, the Stochastic L\"{o}wner 
  Evolution with 
    parameter $\kappa = 6$.  Important and related work by
    Lawler-Schramm-Werner 
    \cite{lsw1, lsw2, lsw3, lsw4, lsw5, lsw6, lsw7} and
    Smirnov-Werner~\cite{sw} has yielded a 
    multitude of results on exponents, conformal invariance and other 
  properties of 
    critical percolation and other two-dimensional processes
(excellent reviews are given in~\cite{lawler, werner2}). 
To extend the work of Schramm and Smirnov, in the spirit of 
Aizenman~\cite{aizenman} 
and Aizenman-Burchard~\cite{ab} (see also~\cite{abnw}), 
it is natural to treat the 
scaling limit for the \emph{set of all contours}, 
as was considered in~\cite{smirnov-long} (see also
Theorem 2.1 of~\cite{lsw5}).
But to our knowledge, no complete description of  
that {\it full\/} scaling limit and its relation to $SLE_6$ has
been proposed, although very interesting ideas do
appear in~\cite{smirnov-long}
(see Theorem~4 and Subsection 3.3 there)
and~\cite{lsw5}
(see Theorem~2.1 there).
In~\cite{cns2}, a certain critical 
dependent percolation model on the triangular (or hexagonal) lattice 
was proved to have the same scaling limit for all of its contours as in 
the independent triangular case, even though the full 
scaling limit itself had not been 
identified. 

In this paper, we present an inductive construction using chordal 
$SLE_6$, which results in a random process (or gas) 
of continuum nonsimple loops in the plane. 
The construction is given 
in Section~\ref{construct} and then a number of features
that we argue are valid for this process 
are presented in Section~\ref{features}. Chief among these features 
is that this process of continuum nonsimple loops 
is indeed the scaling limit (without need for subsequences) of 
the set of all boundary contours for critical site percolation on the 
triangular lattice.
A technical property, which will be used to argue for the
scaling limit feature, is that various ways of organizing
the construction give the same limiting distribution. 
Sketches of the main arguments for the claimed 
features, using~\cite{schramm, smirnov} and other work,
are provided in Section~\ref{limit}. 
A paper by the authors with detailed proofs is in preparation.

  Another important property of the loop process is conformal 
  invariance; 
  we do not discuss that explicitly since it is essentially the same as
  in the conformal invariance results of 
  Lawler-Schramm-Werner~\cite{lsw1, lsw2} 
  and Smirnov~\cite{smirnov} (see also~\cite{lawler, werner2}). 
  We remark that in particular, the distribution of the loop 
  process on the entire plane will be scale and inversion invariant, in 
  addition to being translation and rotation invariant. 

\begin{figure}[!ht]
\begin{center}
\includegraphics[width=8cm]{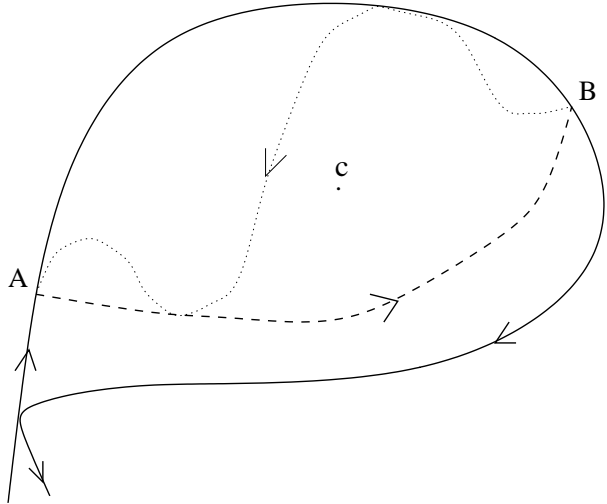}
\caption{Construction of a continuum loop around $c$ in three steps.}
\label{Fig1}
\end{center}
\end{figure}

As a preview of the way in which a single one of our continuum loops is 
constructed, see the (very schematic) Figure~\ref{Fig1}, in which {\it three\/} 
chordal $SLE_6$ processes are used to yield a single loop surrounding a
point $c$ in the plane: The process $\gamma_1$ (solid curve in the
figure), when it first traps $c$ provides a domain 
$D_1$ for the second process $\gamma_2$. A domain $D_2$ for the third 
process $\gamma_3$ (dotted curve) is provided when $\gamma_2$ 
makes an {\it excursion\/} (dashed curve) from $A$ to $B$, two points
on (the ``internal perimeter'' of) $\gamma_1$, and thus traps $c$ 
between itself and (the ``internal perimeter'' of) $\gamma_1$.
The continuum loop consists of 
the excursion segment of $\gamma_2$ from $A$ to $B$ followed by 
$\gamma_3$ from $B$ to $A$.

  We note that in our inductive construction all 
  loops are obtained essentially as in Figure 1, 
  except that the domain $D_1$ may itself be built out of more 
  than one $SLE_6$ process. In the full 
  plane version of our continuum nonsimple loop process, these 
  $SLE_6$ processes 
  (or the excursions into which they
  can be decomposed) are themselves parts of other constructed
  loops and it follows that every new loop touches previous ones
  (e.g., at the points $A$ and $B$ in Figure~\ref{Fig1}).
  This leads to property 4) of 
  Section~\ref{features}, that any two of the continuum nonsimple 
  loops are connected 
  by a ``path'' of touching loops. The analogous lattice result concerns 
  large percolation clusters that almost touch
  (i.e., their boundary contours are separated by only a {\it single\/}
  hexagonal cell) 
  and can be explained in terms of standard
  ``number of arms'' arguments (see~\cite{aizenman}
  and Lemma~5 of~\cite{ksz}). 
It is also related to the high probability that
``fjords'' are of minimal width, a phenomenon observed
numerically by Grossman-Aharony~\cite{ga1, ga2} and explained
by Aizenman-Duplantier-Aharony~\cite{ada},
and which is a key ingredient of our scaling limit claim. 

  In addition to constructing the continuum loops from $SLE_6$ processes, 
  one 
  can also do the reverse --- see property 5) of Section~\ref{features}. 
  We expect 
  that this property, combined with some locality features in the spirit 
  of 
  those already known for $SLE_6$~\cite{lsw1, lawler, werner2}
  should be enough to characterize 
  the full process of continuum nonsimple loops. Other characterizations 
  of the full scaling limit, 
  based on Cardy type crossing formulas~\cite{cardy, cardy2}, have also been 
  proposed --- see, e.g.,~\cite{schramm, smirnov, lsw1, werner}.

\begin{figure}[!ht]
\begin{center}
\includegraphics[width=6cm]{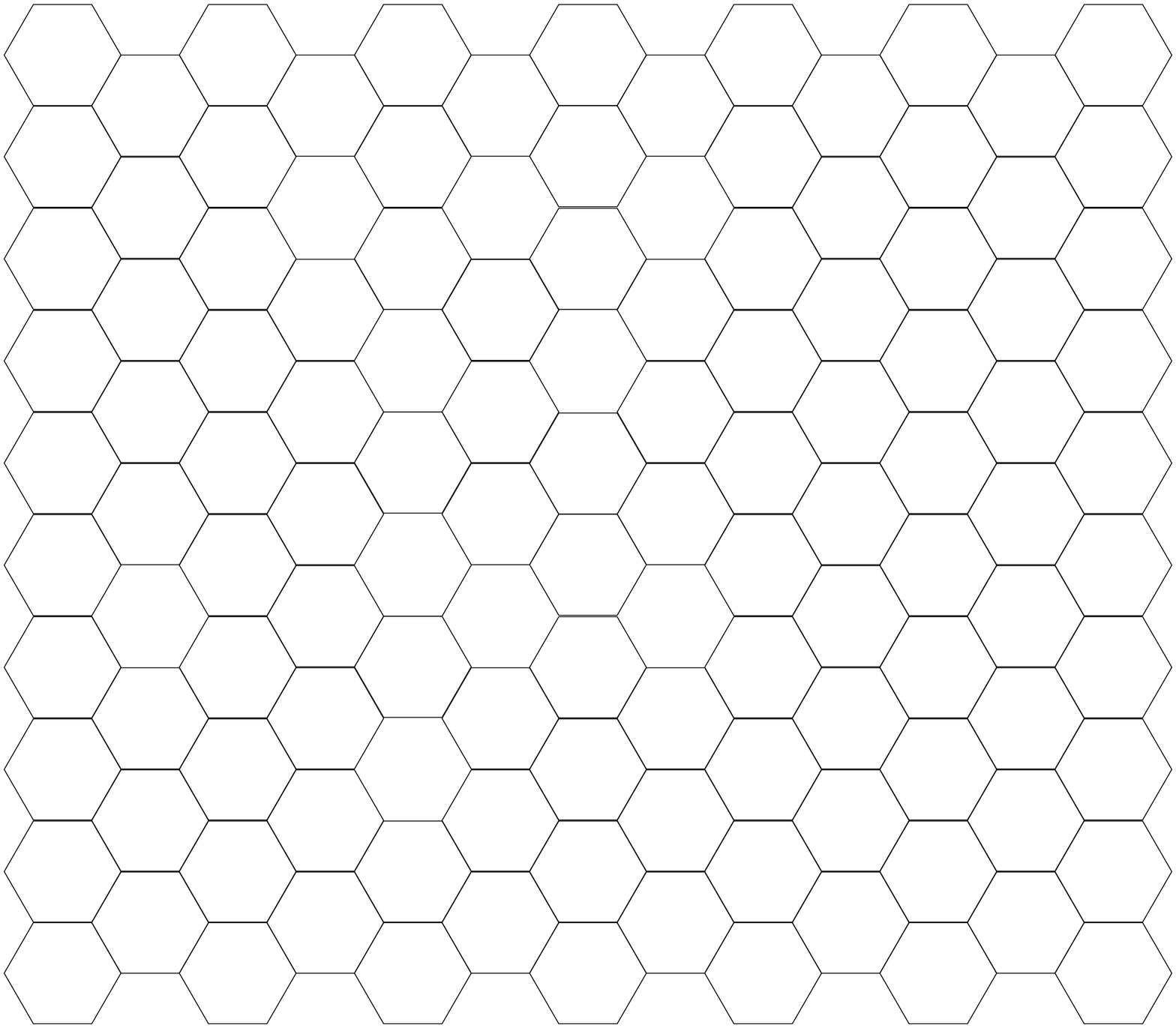}
\caption{Portion of the hexagonal lattice.}
\label{hlattice}
\end{center}
\end{figure}

    The basis for the scaling limit claim,
    presented in Section~\ref{limit}, is 
    a construction for discrete site percolation 
    on the triangular lattice $\mathbb T$, analogous to the construction 
    of the process of continuum nonsimple loops. (We will generally 
  think of 
    the sites of the triangular lattice as the elementary cells of a 
  regular 
    hexagonal lattice $\mathbb H$ embedded in the plane --- see Figure~\ref{hlattice}.) 
    The argument that this discrete 
  construction 
    leads to a proof of the claimed limit is of course itself based 
    on Schramm's~\cite{schramm} and 
    Smirnov's~\cite{smirnov} work on the scaling 
    limit of the percolation ``exploration process,'' which we now briefly 
  review.

    Let $D$ be a bounded simply connected open subset of the plane, with 
    Jordan boundary $\partial D$ (i.e., given by a closed continuous simple curve)
    and two distinct specified points $a,b$ in $\partial D$.
    (Although the restriction to only Jordan regions can probably be
    dispensed with, it is convenient to have it, as we will throughout this
    paper.)
    There is a well-defined stochastic process $\gamma (t) = \gamma_{D,a,b}(t)$
    for $t \in [0,\infty]$ in the closure ${\bar D}$ with $\gamma(0)=a$, 
    $\gamma(\infty)=b$ and H\"{o}lder continuous sample paths that is the 
    trace of Schramm's chordal $SLE_6$, 
    the Stochastic L\"{o}wner Evolution with parameter 
    $\kappa = 6$; this is conventionally defined first on the upper half plane 
    with boundary points $0,\infty$ and then conformally mapped to $D$
    (see~\cite{lsw1}). 

    A major conclusion of Smirnov~\cite{smirnov} is 
    that the scaling limit of a certain ``exploration" process 
    (see Subsection~\ref{exploration} for the definition)
    for critical 
    independent site percolation 
    on the triangular lattice (each site is equally likely to be yellow 
  (minus) 
    or blue (plus)) is 
    the $SLE_6$ process $\gamma$. 
    This is a statement 
    about convergence in distribution, where the topology 
    on sample paths is that of Aizenman-Burchard~\cite{ab}, which 
  uses 
    a supremum norm, but with monotonic reparametrizations of the paths 
  allowed. 
    The exploration process $\gamma^{\delta}$ runs along the edges 
    of the hexagonal lattice that is dual to the triangular lattice of 
  mesh 
  size 
    $\delta$. It 
    basically represents the contour separating blue clusters in $D$ that 
  reach 
    the part $\partial_{a,b} D$ of the boundary $\partial D$ that is 
  traversed 
    when touring $\partial D$ counterclockwise from $a$ to $b$ and 
    yellow clusters 
    in $D$ that reach the other part $\partial_{b,a} D$ of the boundary. 

    The sample paths of $\gamma$ touch both $\partial D$ and themselves 
  many 
    times, but they are noncrossing (and do not touch the same point 
    more than twice or a boundary point even twice). 
    The set $D \setminus \gamma[0,\infty]$ 
    is a countable union of its connected components, which are open and 
    simply connected --- they are Jordan regions like the 
    original region $D$, as will be discussed later 
    in this section (see also~\cite{lsw7, werner2}).
    If $z$ is a deterministic point in $D$, then with probability 
    one, $z$ is not touched by 
    $\gamma$~\cite{rs} and so belongs to a unique one 
    of these, that we denote $D_{a,b}(z)$. 

    There are four kinds of 
    components which may be usefully thought of in terms of how a point 
  $z$ 
    in the interior of the component was first ``trapped'' at some time 
  $t_1$
    by $\gamma[0,t_1]$ perhaps 
    together with either $\partial_{a,b} D$ or 
    $\partial_{b,a} D$: (1) those components whose boundary 
    contains a segment
    of $\partial_{b,a} D$ between 
    two successive visits at $\gamma_0(z)=\gamma(t_0)$ and 
  $\gamma_1(z)=\gamma(t_1)$ 
    to $\partial_{b,a} D$ (where here and below $t_0<t_1$), 
    (2) the analogous components with $\partial_{b,a} D$ 
    replaced by the other part of the boundary $\partial_{a,b} D$, (3) 
    those components formed when $\gamma_0(z)=\gamma(t_0) = 
    \gamma(t_1)=\gamma_1(z)$ with $\gamma$ 
    winding about $z$ in a counterclockwise direction between $t_0$ and 
 $t_1$, 
    and finally (4) the analogous clockwise components.

The boundaries of these components, other than the 
segments of $\partial_{a,b} D$ 
or $\partial_{b,a} D$ in cases (1) and (2), are related to the 
``external perimeter'' of chordal $SLE_6$ 
that was also studied by Smirnov~\cite{smirnov},
Lawler-Schramm-Werner~\cite{lsw7}, and Werner~\cite{werner, werner2}.
    Besides the exploration process itself, there are natural lattice 
  analogues 
    to these components, or more directly relevant for us, lattice 
  analogues 
  to their 
  boundaries
    and to the points $\gamma_0(z),\gamma_1(z)$ on their boundaries. 
    We argue that it should follow from the work of Smirnov 
    combined with other percolation arguments 
    (see Subsection~\ref{convergence})
    that for any finitely (or countably) many deterministic points 
  $z_1,z_2,\dots$ 
    in $D$, 
    the joint distribution of 
  the corresponding boundaries and points 
  converges to the distribution of the 
  continuum $SLE_6$ objects. 
  This convergence also shows that the boundaries of these
  regions are Jordan curves (see~\cite{ada} and also~\cite{lsw7, werner2}).

  To obtain these 
  new lattice analogues, one ``fattens'' the exploration process from 
  being a path along the dual lattice (i.e., along the edges
  of $\delta {\mathbb H}$) to include all the 
  blue and yellow sites that touch that path (i.e., the hexagons that 
  have 
  actually been explored while constructing the path $\gamma^{\delta}$). 
  Then one considers the 
  connected (in the lattice sense) components of the difference between 
  the set of all the sites in $D$ and the set of all sites in the fattened 
  exploration path. 
  For more details, see Section~\ref{limit}.

  \section{Construction of the Continuum Nonsimple Loops} 
  \label{construct}

In defining our process, we will freely switch between the real
plane ${\mathbb R}^2$, with $x$ and $y$ coordinates, and the complex
plane $\mathbb C$, according to convenience.
The basic ingredient in the algorithmic construction consists of a chordal 
$SLE_6$ path between two points $a$ and $b$ of the boundary $\partial D$ of 
a given simply connected domain $D \subset {\mathbb C}$. 
The domains we will encounter in the construction are bounded open sets $D$ 
whose boundaries $\partial D$ are Jordan curves. 
This ensures that the unit disc $\mathbb U = \{ z \in {\mathbb C} : |z| < 1 \}$ 
can be mapped onto $D$ via a conformal transformation that can be extended 
continuously to the boundary.

As we will explain soon, sometimes the two boundary points are
``naturally'' determined as a product of the construction itself,
and sometimes they are given as an input to the construction.
In the second case, there are various procedures which would yield
the ``correct'' distribution for the resulting continuum nonsimple
loop process; one possibility is as follows.
Given a domain $D$, choose $a$ and $b$
so that, of all points in $\partial D$, they have maximal $x$-distance
or maximal $y$-distance, whichever is greater.
A crucial aspect of this procedure, as discussed in Subsec.~\ref{finding}
below, is that there is a bounded away from zero probability that the
resulting subdomains $D_{a,b}(z)$ have maximal $x$-distance (or else
maximal $y$-distance) shrunk by a bounded away from one factor compared
to the domain $D$. Another aspect, implicit in 
Subsecs.~\ref{full} and \ref{convergence}, is that the corresponding
$(a,b)$'s of our discrete construction converge in distribution (in the
scaling limit) to those of the continuum construction. This latter
aspect and also the well-definedness of the above procedure would
be resolved by showing that, in the context of our continuum construction, 
the choices of $(a,b)$ are unique with probability one (for the starting
domain, the unit disc ${\mathbb U}$, we take $(a,b)=(-i,i)$). We believe
that this is so and proceed under that assumption, but in any case, 
the issue can be avoided by doing
a randomized version of the above procedure in which $(a,b)$ are chosen
to be ``fairly close'' to having the maximal $x$-distance or $y$-distance.

To start our construction, we take the unit disc ${\mathbb U} = {\mathbb U}_1$
(later, to take a thermodynamic limit and extend the loop process 
to the entire plane, the unit disc will be replaced by a 
growing sequence of large discs, ${\mathbb U}_R$) and
begin by ``running'' chordal $SLE_6$ inside ${\mathbb U}$ 
from $a=-i$ to $b=i$. 

  The resulting path $\gamma_{{\mathbb U},-i,i}$ (the trace 
  of chordal 
  $SLE_6$) touches itself and $\partial {\mathbb U}$ (infinitely) many 
  times. 
  The set ${\mathbb U} \setminus \gamma_{{\mathbb U},-i,i}[0,\infty]$ 
  is a 
  countable union of its connected components, which are open and simply 
  connected. 
  They can be of four different types, as explained in the introduction. 

  To conclude the first step (in this version of the
  construction), we consider all domains of type (1), 
  corresponding to excursions of the $SLE_6$ path from the left portion 
  of the boundary of the unit disc (the one from $i$ to $-i$
  counterclockwise). 
  For each such domain $D$, the points $a$ and $b$ on its boundary are 
  chosen 
  to be respectively those points where the excursion ends and where 
  it begins, that is, if $z \in D$, we set $a=\gamma_1(z)$ and 
  $b=\gamma_0(z)$ (in the notation of Section~\ref{intro}). 
  We then run chordal $SLE_6$ from $a$ to $b$. 
  The loop obtained by pasting together the excursion 
  from $b$ to $a$ followed by the new $SLE_6$ 
  path from $a$ to $b$ is one of our continuum loops. 
  At the end of the first step, then, the procedure has generated 
  countably 
  many loops that touch the left side of the original boundary (the 
  portion 
  $\partial_{i,-i} {\mathbb U}$ of the boundary of $\mathbb U$); 
  each of these loops touches the left side of the original
  boundary but may or may not touch the right side.

The last part of the first step also produces new domains, corresponding 
to the connected components of $D \setminus \gamma_{D,a,b}[0,\infty]$ for 
all domains $D$ of type (1). 
Each one of these components, together with all the domains of type (2),
(3) and (4) previously generated, is to be used in the next step of the 
construction, playing the role of the unit disc. 
For each one of these domains, we choose the ``new $a$''
and ``new $b$'' on the boundary as explained before, and then continue
with the construction.
Note that the ``new $a$'' and ``new $b$'' are chosen according to the
rule explained at the beginning of this section also for domains of
type (2), even though they are generated by excursions like the domains
of type (1).

  This iterative procedure produces at each step a countable set of loops. 
  The limiting object, corresponding to the collection of all such loops, 
  is our basic process. 
  (Technically speaking, we should include 
  also trivial loops fixed at each $z$ in $\mathbb C \cup \{\infty\}$ so 
  that the collection of loops is closed in an appropriate sense~\cite{ab}.) 
  Some of its properties will be given in the next section. 

  As explained, the construction carries on iteratively and can be 
  performed 
  simultaneously on all the domains that are generated at each step. 
  We wish to emphasize, though, that the obvious monotonicity of the 
  procedure, where at each step new paths are added
  independently in different domains, and new domains are 
  formed from the existing ones, implies that any other choice of the 
  order 
  in which the domains are used would give the same result (i.e., produce 
  the same limiting distribution), provided that every 
  domain that is formed 
  during the construction is eventually used. 
In Section~\ref{limit}, when arguing that the lattice scaling 
limit coincides with our continuum nonsimple loop process, it will be
convenient to utilize a different 
procedure in which each step involves only a single
$SLE_6$ for a single domain. 
This will be done with the help of a deterministic set of points $\cal P$
that are dense in $\mathbb C$ and are endowed with a deterministic order.
The domains will then be used one at a time, with domains containing
higher ranked points of $\cal P$ having a higher priority for 
order of being used. 

  In Section~\ref{limit}, arguments will be given as to why the process of 
  loops 
  we have just constructed is the scaling limit, as 
  $\delta \to 0$, of the set of all cluster boundary 
  contours for critical percolation on the portion of the 
  triangular lattice of mesh size $\delta$ sitting within the disc 
  ${\mathbb U}_1$ 
  of radius $1$, and with blue (plus) boundary conditions imposed. Of 
  course, 
  essentially the same construction and scaling limit results 
  can be done on the disc ${\mathbb U}_R$ of radius $R$. 
  It is not hard to then verify that the limit in distribution 
  of the loop process exists as $R \to \infty$ 
  and that this represents the scaling limit of the set of all cluster 
  boundary contours 
  in the entire plane, with no boundary conditions needed. It is this 
  process in 
  the entire plane that we will consider in the next section of the paper 
  dealing 
  with properties of the loop process. 

  \section{Features of the Continuum Nonsimple Loop Process} 
  \label{features} 

  In this section we present a number of features that we 
  argue will 
  be valid for our process of continuum nonsimple loops in the plane. 
  Some of them are direct
  consequences of the continuum algorithmic 
  construction, while others become clear only in light of the analogous 
  construction for discrete percolation, which will be presented in 
  the next section.
%
%
%

The first feature is the scaling 
limit property --- which is
used to derive the other properties.
A sketch of the derivation of the scaling limit and other properties 
is given in the next section of the paper. 
The scaling limit property~1) 
is a distributional statement;
properties~2)--4) all involve  statements that are valid with probability one;
property 5) is a bit of a hybrid.

  \bigskip 
  \noindent 
  {\bf 1)} The continuum nonsimple loop 
  process is the scaling limit of the set of all boundary 
  contours for critical site percolation on the triangular lattice.

  \bigskip 
  \noindent 
  {\bf 2)} This process is a random collection
  of noncrossing, continuous loops 
  on the plane. 
  The loops touch themselves and each other many times, but no three 
  or more loops can come together at the same point, and a single loop 
  cannot touch the same point more than twice.

  \bigskip 
  \noindent 
  {\bf 3)} Any deterministic point of the plane is 
  surrounded by an infinite family 
  of nested loops with diameters going to both zero 
  and infinity; any annulus about 
  that point with inner radius $r_1 > 0$ and outer radius $r_2 < \infty$ 
  contains only a finite number of loops. 
  Consequently, any 
  two distinct deterministic points of the plane are separated by loops 
  winding 
  around them. 

  \bigskip 
  \noindent 
  {\bf 4)} Any two loops are connected by a finite ``path'' of touching loops. 

  \bigskip 
  \noindent 
  {\bf 5)} For a (deterministic) Jordan region $D$ with two boundary points $a$ 
  and $b$, 
  a process distributed as chordal $SLE_6$ from $a$ to $b$ can be 
  constructed 
  starting from the continuum nonsimple loops (in the whole plane) by 
  doing a 
  continuum analogue of what is done on the lattice to piece together 
  cluster boundary 
  segments 
  to give the lattice percolation exploration process (see below). 

  \bigskip 

  We conclude this section of the paper with a more detailed 
  explanation of the construction just mentioned in property 5). 
  To do the construction, 
  it is useful to first convert all the loops into directed ones. 
  There is one binary choice to be made: any one loop can be given either 
  the 
  clockwise or counterclockwise direction and then all other loops are 
  automatically 
  determined (via the natural nesting tree structure of the set of all the 
  loops) 
  by requiring that the set of all loops (ordered by nesting) about any 
  deterministic 
  point alternate in direction. Back on the percolation lattice the two 
  choices 
  correspond to either having yellow just to the left of the 
  directed path and blue just to 
  the right or vice versa; the two choices are also of course related by a 
  global color flip. 

  For convenience, 
  let us suppose that $a$ is at the bottom and $b$ is at the top
  of $D$ so that the boundary 
  is divided into a left and right part by these two points. 
  The desired path from $a$ to $b$ is then put together using 
  {\it all\/} of the 
  following directed segments of the loops 
  (most of the analogous types of segments for the lattice
  exploration process may be seen in Figure~\ref{Fig3}):
 (i) ``excursions'' from the left part 
  to the right part of the boundary (they touch each of the 
  two boundary parts at exactly 
  one point), (ii) ``excursions'' from right to left, 
  (iii) excursions from the left part to itself (they touch that part of 
  the 
  boundary at exactly two points) which do not touch the right part and 
  which are 
  maximal in that there is not another such excursion between them and the 
  right part, 
  and (iv) the analogous excursions from the right part to itself. 

  There are countably many
  excursions of types (i) and (ii) which are ordered from lower to 
  higher and alternate 
  between types (i) and (ii). The segment of the right boundary between 
  where an 
  excursion of type (i) ends and the next excursion of type (ii) begins 
  supports 
  countably many excursions of type (iv) which are also ordered from 
  lower to 
  higher. These may be all pieced together (they don't quite touch 
  so a limit is needed) in order that a continuous path 
  connecting 
  the type (i) to the type (ii) excursion is obtained. Using such 
  connecting paths 
  on the right and the analogous paths on the left that connect the end 
  of a type (ii) to the beginning of the next type (i), one can connect 
  all the type (i) and (ii) excursions in order and obtain finally the 
  desired 
  path from $a$ to $b$. 

  \section{The Continuum Nonsimple Loop Process as Scaling Limit} 
  \label{limit} 

  In this section we will introduce a discrete inductive 
  construction which is analogous to the continuum construction 
  given in Section~\ref{construct}. 
  Our interest in the discrete construction comes from the claim 
  that the continuum one is its scaling limit. 
  This requires comparing the two constructions. 
  In order to do so, we first reorganize the continuum one 
  and introduce some notation. 

  \subsection{Priority-Ordered Continuum Construction} 
  \label{priority} 

  We want to arrange the continuum construction in such a 
  way that each step corresponds to a single new $SLE_6$ path. 
  To do that, we need to order the domains present at the beginning 
  of each 
  stage (which is the term we use for a group 
  of successive single steps), 
 so as to choose which ones to use in the steps of that stage. 
  The domains are the connected components that the original domain 
  is broken up into by {\it all\/} the $SLE_6$ paths constructed 
  up to the beginning of the new stage. 
  The ordering will be done with the help of the deterministic 
  ordered set 
  of points $\cal P$, dense in $\mathbb C$, introduced in 
  Section~\ref{construct}. 

  The first step and stage consists of an $SLE_6$ path 
  $\gamma_1 = \gamma_{{\mathbb U},-i,i}$ 
  inside $\mathbb U$ from $-i$ to $i$ which, as explained in 
  Section~\ref{construct}, produces many domains which are the 
  connected components of the set $\mathbb U \setminus 
  \gamma_1 [0,\infty]$. 
  These domains can be priority-ordered using points in $\cal P$, 
  according to the rank of the highest ranking point of 
  $\cal P$ that each contains. The priority orders of domains change 
  as the construction proceeds. 

The second stage of the construction consists of two $SLE_6$ paths, 
$\gamma_2$ and $\gamma_3$, that are produced in the two domains 
with highest priority at the end of the first stage,
the priority being determined using the points
of $\cal P$ and the starting and ending points for domains 
that are not of type~1) being chosen as explained in Section~\ref{construct}. 

  In general, for the $k$th stage of the construction, $k$ 
  $SLE_6$ paths are produced in those $k$ domains present at the 
  end of the last stage 
  with highest priority, again using the points of $\cal P$ 
  for ranking the domains.
  This way of organizing the construction does not affect 
  the final result, as discussed in Section~\ref{construct},
  and has the advantage that to each step corresponds a single 
  $SLE_6$ path, with the $SLE_6$ paths ordered. 

\subsection{Discrete Exploration and Loop Construction} 
\label{exploration} 

We will organize the discrete construction, which we will 
present soon, in the same way. 
Before doing that, though, we briefly introduce its 
key ingredient --- the discrete exploration 
process for a general simply connected set $D^{\delta}$ of hexagons. 

To begin, we denote by $\partial D^{\delta}$ the edge boundary.
For two points, $a, b$ in $\partial D^{\delta}$ suitably chosen
at the vertices of two hexagons, the usual exploration
process~\cite{smirnov} (see also~\cite{lawler, werner2, werner})
with $\pm$ boundary conditions (i.e.,
blue hexagons just outside the counterclockwise portion 
$\partial_{a,b} D^{\delta}$ of $\partial D^{\delta}$ from $a$ to $b$
and yellow hexagons just outside the other portion $\partial_{b,a}
D^{\delta}$)
can be described as a sort of  self-avoiding random walk on the edges
of the hexagons contained in $D^{\delta}$ that moves left
(with respect to the current direction of exploration)
when a blue hexagon is encountered and right when a yellow one is
encountered.

\begin{figure}[!ht]
\begin{center}
\includegraphics[width=8cm]{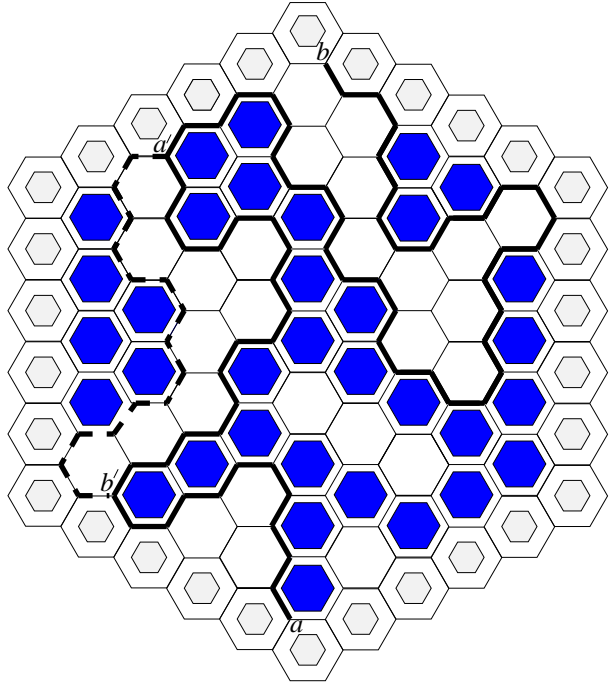}
\caption{Construction of the outer contour of a cluster
of yellow/minus (white in the figure) hexagons in two steps
by means of a first exploration from the vertex $a$
to $b$ (heavy line), followed by a second one from
$a'$ to $b'$ (heavy broken line).
The outer layer of hexagons does not belong to the domain
where the explorations are carried out, but represents its
monochromatic blue/plus external site boundary.}
\label{Fig3}
\end{center}
\end{figure}

We use this rule for $\pm$ boundary conditions, and \emph{also}
for $+$ (blue) boundary conditions, proceeding at the boundary
\emph{as if} we had $\pm$ boundary conditions
(see Figure~\ref{Fig3}).
For $\mp$ and $-$ (yellow) boundary conditions, we use the ``opposite''
(with respect to color) rule.
We remark that although the exploration process itself changes under a 
color flip of the boundary conditions, its distribution is
color-blind. 

The interpretation of the exploration process depends
on whether the boundary condition is monochromatic or not.
Let $\Delta D^{\delta}$ be the external (outer) site
boundary of $D^{\delta}$, with $\Delta_{a,b} D^{\delta}$
and $\Delta_{b,a} D^{\delta}$ representing the portions
next to $\partial_{a,b} D^{\delta}$ and $\partial_{b,a} D^{\delta}$
respectively.

  \begin{itemize} 
  \item For regions with $\pm$ 
  (respectively, $\mp$) boundary conditions, the exploration 
  path represents the contour separating the blue 
  (respectively, yellow) cluster that 
  contains $\Delta_{a,b} D^{\delta}$ from the yellow 
  (respectively, blue) cluster that 
  contains $\Delta_{b,a} D^{\delta}$. 
  \item For regions with monochromatic blue 
  (respectively, yellow) boundary conditions, the 
  exploration path represents portions of the outer boundary contours of 
  yellow (respectively, blue) clusters touching $\partial_{b,a} D^{\delta}$ 
  and adjacent to blue 
  (respectively, yellow) hexagons that are the starting point of 
  a blue (respectively, yellow) path (possibly an empty path) 
  that reaches 
  $\partial_{a,b} D^{\delta}$, 
  pasted together using 
  portions of $\partial_{b,a} D^{\delta}$. 
  \end{itemize} 


Next, we show how to get the complete outer 
contour of a monochromatic (say, yellow) 
cluster by twice using the exploration process described above
(see Figure~\ref{Fig3}).
Consider a large simply connected domain $D^{\delta}$ surrounded
by blue hexagons, which we can identify with $\Delta D^{\delta}$.
$D^{\delta}$ will contain many clusters of both colors in its interior. 
We pick two suitably chosen points $a, b \in \partial D^{\delta}$
and perform the exploration from $a$ to $b$. 

While performing the exploration process, we discover the 
color of the hexagons that touch the exploration path. 
We want to keep track of that information. 
As a result, at the end of the exploration process we have 
three ``paths'': the exploration path $\gamma^{\delta}$ along 
the edges of the hexagonal lattice, and respectively the 
``paths''  $\Gamma^{\delta}_{Y}$ and $\Gamma^{\delta}_{B}$ 
along the (respectively, yellow or blue) 
sites of the triangular lattice that touch it (i.e., those 
hexagons that have at least one edge belonging to the 
exploration path). 
The latter lattice ``paths'' are not in general simple,
as they can form loops and have dangling ends. 

The set $D^{\delta} \setminus \{ \Gamma^{\delta}_{Y} \cup 
\Gamma^{\delta}_{B} \}$  is the union of its connected 
components (in the lattice sense), 
which are simply connected.
There are four types of components which may be usefully thought 
of in terms of their external site boundaries: (1) those components 
whose  site boundary contains both sites in $\Gamma^{\delta}_{Y}$ 
and  $\Delta_{b,a} D^{\delta}$, (2) the analogous components 
with $\Delta_{b,a} D^{\delta}$ replaced by  $\Delta_{a,b} D^{\delta}$ 
and $\Gamma^{\delta}_{Y}$  by $\Gamma^{\delta}_{B}$, 
(3) those components whose site  boundary only contains sites 
in $\Gamma^{\delta}_{Y}$, and  finally (4) the analogous components 
 with $\Gamma^{\delta}_{Y}$  replaced by $\Gamma^{\delta}_{B}$. 

If we now take a region of type (1), there are natural 
starting and ending points 
(where the excursion
that produces that region respectively
ends and starts; e.g., $a', b'$ in Figure~\ref{Fig3}) 
for an exploration process within it.
Performing such an exploration process inside the specified 
domain of type (1) and pasting the new exploration path together 
with the portion of a previous exploration path corresponding
to the excursion that produced that domain of type (1)
will generate a loop along the edges of the hexagonal lattice. 
The loop is the outer contour of a yellow cluster 
that touches $\partial_{b,a} D^{\delta}$ and is adjacent 
(on its ``right'')
to blue hexagons, each of which is the starting point of a blue 
path to $\partial_{a,b} D^{\delta}$. 

Analogous exploration processes in the other regions of type (1) 
produce similar loops on the edges of $\delta \mathbb H$ that are also
boundary contours.
In fact, every domain with $\pm$ (or $\mp$) boundary conditions
obtained during the discrete algorithmic construction that we are about 
to present will contain an exploration path which, pasted together 
with the appropriate part of a previous exploration path, 
provides the complete outer boundary contour of a monochromatic cluster. 

\subsection{Full Discrete Construction} 
\label{full} 

We now give the algorithmic construction for discrete percolation 
which is the analogue of the continuum one. 
Each step of the construction is a single percolation exploration 
process; the order of successive steps is organized as in the continuum 
construction detailed at the beginning of this section. 
We start with the set $D^{\delta}_0$ of hexagons that are 
contained in the unit disc $\mathbb U$
and will make use of the deterministic countable 
ordered set $\cal P$ of points 
dense in $\mathbb C$ that was introduced in Section~\ref{construct}. 

The first step consists of an exploration process inside 
$D^{\delta}_0$. 
For this, we need to select two points $a$ and $b$ in 
$\partial D^{\delta}_0$.
We choose for $a$ some vertex close to $-i$, and for $b$ one
close to $i$.
The first exploration produces a path $\gamma^{\delta}_1$ 
and, for $\delta$ small, many new domains of all four types. 
These domains are ordered with the help of points in $\cal P$ 
as in the continuum case, and that order is used, at each 
stage of the process, to determine the next group of exploration processes. 
So, for the second stage of the construction, two domains are 
chosen and explored, and so on. 
With this choice, the exploration processes and paths are 
naturally ordered: $\gamma^{\delta}_1, \gamma^{\delta}_2, \ldots$ . 

Each exploration process of course requires choosing a starting 
and ending point, which is done mimicking what is done in the
continuum case (with some adjustments due to the discrete nature
of the lattice).
For domains of type (1), with $\pm$ or $\mp$ boundary conditions,
the choice is the natural one, explained before, which produces 
a loop using the edges of $\delta \mathbb H$.
For a domain $D^{\delta}_k$ (used at the $k$th step) of type 
other than (1), and therefore with monochromatic boundary
conditions, two vertices are chosen that are close to the
two points of $\partial D^{\delta}_k$ selected according
to the rule given in Section~\ref{construct}.

The procedure continues iteratively, with regions of type (2), 
(3) and (4), which have monochromatic boundaries, playing 
the role played in the first step by $D_0^{\delta}$.
As the construction continues, new loops along the edges 
of the hexagonal lattice are formed which correspond to the outer 
boundary contours of constant sign (monochromatic) clusters. 

\subsection{Ingredients for Convergence} 
\label{convergence} 
By comparing the discrete and continuum version of the algorithmic
construction, and using repeated applications of Smirnov's
work~\cite{smirnov}, we will argue
that for any fixed $k$, the first $k$ steps of the
discrete construction converge (jointly, in distribution)
to the first $k$ steps of the continuum construction, as $\delta \to 0$.
This claim is an extension of the discussion near the end
of Section~\ref{intro} about convergence in distribution of certain lattice 
boundaries and points to their continuum analogues.
We note that one complication is due to the fact
that the boundaries of the domains where the exploration processes
are performed are not deterministic, but are
themselves obtained using exploration processes.
Some continuity arguments are therefore needed.

%
%

\subsubsection{Matching Continuous and Discrete Domains and Loops}
\label{matching}

A key ingredient is the observation that the probability
of ``fjords'' of width larger than the minimal one goes
to zero in the scaling limit~\cite{ada}.
This ensures that the domains and loops
generated at various steps of
the continuum construction are the limits of
corresponding domains and loops produced
in the discrete one, so that, e.g., one can identify, with probability
going to one as $\delta \to 0$, the domain containing a point
$c$ at a given step of the continuum construction with the
domain containing $c$ at the equivalent step of the discrete one.

\subsubsection{Finding Large Contours in $O(1)$ Steps} 
\label{finding} 

The discrete algorithm will reach and discover 
all the boundary contours inside $\mathbb U$;
moreover we argue that the 
number of steps $K_{\varepsilon} (\delta)$ needed for the discrete 
algorithm to recover all contours in ${\mathbb U}$ 
of diameter larger than a given 
$\varepsilon > 0$ is bounded (in probability) as $\delta \to 0$. 

This uses the observation that the discrete algorithm cannot ``skip''
a contour and move to explore the domain inside it and the fact that
the maximum diameter of the domains present inside
$\mathbb U$ after $k$ steps of the discrete algorithm tends to zero
in probability as $k \to \infty$, $\delta \to 0$.
To understand the last fact, first of all notice that
the construction cannot produce ``too many'' distinct domains
of diameter greater than $\varepsilon$, or else there would be
too many disjoint ``macroscopic'' monochromatic paths
(the site boundaries of those domains) in
$\delta {\mathbb T} \cap {\mathbb U}$ to satisfy
the multiple crossing probability bounds of~\cite{ab}.
Consider now a domain $D^{\delta}$ with points $a$ and $b$
on the boundary $\partial D^{\delta}$ chosen because they have,
among all points in $\partial D^{\delta}$, maximal $x$-distance.
Then standard percolation arguments~\cite{russo, sewe} ensure
that, with bounded away from zero probability, the maximal
$x$-distance between points on the boundary of each of the
components that $D^{\delta}$ is split up into by effect of the exploration
process is smaller than, say, two thirds of the $x$-distance between
$a$ and $b$.
(Notice also that each newly formed domain is ``unexplored territory''
on which no information is available before the exploration process
inside it begins.)

The proof of property 1) of Section~\ref{features} is
completed by first letting $\varepsilon \to 0$ and
then by taking the thermodynamic limit 
(to obtain a loop process in the entire plane, as discussed at the 
end of Section~\ref{construct}). 
  
\subsection{Properties of the Continuum Loop Process} 
\label{properties} 

We now turn to brief sketches of the derivations
of the other properties presented in Section~\ref{features}.


\bigskip
\noindent 
{\bf 2)} The noncrossing property of contours is preserved in the 
scaling limit, and the fact that they touch themselves and 
each other follows fairly directly from the continuum
construction (see the discussion below about property~4)). 
The properties that no three or more loops can come together 
at the same point and a single loop cannot touch the same point 
more than twice follow from standard ``number of arms'' arguments
(see~\cite{aizenman} and Lemma~5 of~\cite{ksz}). 

\bigskip 
\noindent 
{\bf 3)} Both the fact that any deterministic point of the
plane is surrounded by infinitely many loops and the claim about the inner
and outer radii of loops surrounding a given point follow from property~1)
combined with standard percolation arguments~\cite{russo, sewe} 
(see also Lemma~3 of~\cite{ksz}).

\bigskip 
\noindent 
{\bf 4)} This property follows fairly directly from the continuum
construction, as discussed in Sec.~\ref{intro}. 
As explained in the introduction, the analogous lattice result
concerns large clusters of the same sign that almost touch and
the existence of ``macroscopic fjords'' only of minimal
width (see~\cite{ga1, ga2, ada}).
For example, the existence of a long double monochromatic
layer of hexagons separating two large clusters of the same
color would give rise to six disjoint ``macroscopic''
paths of hexagons not all of the same color
which start within a ``microscopic'' distance of each other.
The probability of this happening goes to zero as $\delta \to 0$.

\bigskip 
\noindent 
{\bf 5)} This property is proved by noting that the usual lattice exploration 
process can be realized as a discrete version of the 
continuum exploration 
procedure outlined at the end of Section~\ref{features}. 
By~\cite{smirnov}, it is enough to show that the lattice version 
converges to the continuum one.

\bigskip 

\noindent {\bf Acknowledgments.} We benefited greatly from 
attending the May, 2002 lectures on 
Conformally Invariant Processes given by Greg Lawler at the 
ICTP School of Probability Theory in Trieste 
and from conversations with him about SLE 
and related issues. We thank Bill Faris for useful suggestions 
about the presentation of our results. F.~C.~thanks the Courant 
Institute 
for its kind hospitality during the period when this work was carried 
out.

\end{document}